\documentclass[11pt]{article}
\usepackage{amsmath}
\usepackage{amsmath,amssymb,amsfonts,amsthm,fancyhdr}
\usepackage{epsfig,graphicx,picins,picinpar,subfigure}
\usepackage{pstricks}
\usepackage{fancyvrb}
\usepackage[numbers,sort&compress]{natbib}

\begin{document}

\title{\textbf{The problem of the least prime number in an arithmetic progression and its applications to Goldbach's conjecture}}
\author{SHAOHUA ZHANG}
\date{{\small\emph{School of Mathematics, Shandong University,
Jinan,  Shandong, 250100, PRC\\
E-mail address: shaohuazhang@mail.sdu.edu.cn}}} \maketitle

\vspace{3mm}\textbf{Abstract:} The problem of the least prime number
in an arithmetic progression is one of the most important topics in
Number Theory. In [11], we are the first to study the relations
between this problem and Goldbach's conjecture. In this paper, we
further consider its applications to Goldbach's conjecture and
refine the result in [11]. Moreover, we also try to generalize the
problem of the least prime number in an arithmetic progression and
give an analogy of Goldbach's conjecture.

\vspace{3mm}\textbf{Keywords:} least prime number, arithmetic
progression,  Kanold's hypothesis, Chowla's hypothesis, Goldbach's
conjecture, Dirichlet's theorem, Dickson's conjecture, Chinese
Remainder Theorem, prime map, prime number,  prime point

\vspace{3mm}\textbf{2000 MR  Subject Classification:} 11A41, 11A99,
11B25, 11P32


\section{Introduction}
\setcounter{section}{1}\setcounter{equation}{0}
Let  $k,l$ denote positive integers with $(k,l)=1$ and $1\leq l\leq
k-1$. Denote by $p(k,l)$ the least prime $p\equiv l(\mod k)$. Let
$p(k)$ be the maximum value of $p(k,l)$ for all $l$ with $(k,l)=1$
and $1\leq l\leq k-1$. In 1944, Linnik [1] proved that $p(k)<k^L$,
$L$ is an absolute constant which is now called Linnik's constant.
In 1957, Pan [2] claimed $L\leq10000$. In 1958, he [3] was the first
to prove that $L\leq 5448$. In 1992, Heath-Brown [4] proved $p(k)\ll
k^{5.5}$.  Recently,  Xylouris [5] improved this result to $p(k)\ll
k^{5.2}$. In 1989, Bombieri, Friedlander and Iwaniec [6] proved
$L\leq 2$ for almost all integers. Kanold [7, 8] ( also
independently made by Schinzel and Sierpi\'{n}ski [9]) conjectured
that $p(k)<k^2$ for every positive integer $k>1$. In [4],
Heath-Brown [4] proved $p(k)<(\varphi (k)\log k)^2$ assuming the
Generalized Riemann Hypothesis.   Chowla [10] has observed that
$p(k)\ll k^{2+\epsilon}$ for every $\epsilon>0$ assuming the
Generalized Riemann Hypothesis. He further conjectured $p(k)\ll
k^{1+\epsilon}$ for every $\epsilon>0$. Thus, we have the following
weakened form of  Chowla's hypothesis:

\vspace{3mm}\noindent{\bf  Conjecture 1:~~}%
For any positive real number $0<\varepsilon<\frac{1}{2}$, there is a
positive constant $C_1$ depending on $\varepsilon$ such that for
every sufficiently large positive integer $k>C_1$, $p(k)<
k^{2-\varepsilon}$.

\vspace{3mm}From the aforementioned rich achievements and
advancements, one see that the problem of the least prime number in
an arithmetic progression is very interesting. It is one of the most
important topics in Number Theory. In 2008, we [11] found that this
problem closely ties up Goldbach's conjecture. In this paper, we try
to refine the result in [11]. Moreover, we also try to generalize
the problem of the least prime number in an arithmetic progression
and give an analogy of Goldbach's conjecture.

\section{Main results}
\vspace{3mm}\noindent{\bf  Lemma 1:~~}%
For any integer $n>6$, there must be two distinct odd primes $p,q$
such that $\gcd(pq,n)=1$ and $p<n,q<n$.

\vspace{3mm}\noindent{\bf  Proof of Lemma 1:~~}%
By the refined Bertrand-Chebyshev theorem which states that there
exists at least two distinct primes in the interval $(m, 2m)$ when
$m = 4$ or $m > 5$, it is easy to prove Lemma 1 holds since any
prime in the interval $(\frac{n}{2}, n)$ is coprime to $n$.

\vspace{3mm} Goldbach's famous conjecture states that every even
integer $2n\geq 4$ is the sum of two primes. Due to it is trivial
that it is true for infinitely many even integers:  $2p=p+p$ (for
every prime $p$), we give Goldbach's conjecture a slightly different
expression that every even integer $2n\geq 8$ is the sum of two
distinct primes. Thus, by Lemma 1, we get a necessary condition of
Goldbach's conjecture as follows.

\vspace{3mm}\noindent{\bf  Conjecture 2:~~}%
For integer $n>6$, there exists a natural number $r$  such that
$2n-p_r$ is coprime to each of
$2n-p_1,...,2n-p_{r-1},2n-p_{r+1},...,2n-p_k$, where
$p_1,...,p_{r-1},p_r, p_{r+1},...,p_k$ are all old primes smaller
than $n$, $p_r$ satisfies $\gcd(p_r,n)=1$ and $1\leq r\leq k=\pi
(n-1)-1$, where $\pi(x)$ is the prime counting function giving the
number of primes less than or equal to a given number $x$.

\vspace{3mm}\noindent{\bf  Theorem  1:~~}%
Conjecture 1 and Conjecture 2 imply that every sufficiently large
even integer may be written as the sum of two distinct primes.

\vspace{3mm}\noindent{\bf  Lemma 2:~~}%
Denote the least prime coprime to $m$ by $q(m)$ for any positive
integer $m$. For every integer $k\geq 1$, there is a positive
integer $C_2$ depending on $k$ such that for every integer $m\geq
C_2$, we have $(q(m))^k<m$.

\vspace{3mm}\noindent{\bf  Proof of Lemma 2:~~}%
If $k=1$, the proof is straightforward. Let's consider the case
$k>1$. By P\'{o}sa's result [12] or Prime Numbers Theorem, there is
a positive integer $n_k$ such that $p_{n+1}^k<p_1p_2\cdots p_n$ for
all $n\geq n_k$, where $p_i$ is the $i$th prime. Let
$C_2=p_1p_2\cdots p_{n_k}$. We claim that for every integer $m\geq
C_2$, we have $(q(m))^k<m$. Write $p_1p_2\cdots p_r\leq
m<p_1p_2\cdots p_{r+1}$. Since $m\geq C_2=p_1p_2\cdots p_{n_k}$,
hence $r\geq n_k$. So, $p_{r+1}^k<p_1p_2\cdots p_r$. If $q(m)\leq
p_{r+1}$, then $(q(m))^k\leq p_{r+1}^k<p_1p_2\cdots p_r\leq m$ and
Lemma 2 holds. If $q(m)> p_{r+1}$, then $m$ is divisible by
$p_1p_2\cdots p_{r+1}$ because $q(m)$ is the least prime coprime to
$m$. Therefore, $m\geq p_1p_2\cdots p_{r+1}$. It is a contradiction
by our assumption on $p_1p_2\cdots p_r\leq m<p_1p_2\cdots p_{r+1}$.
So Lemma 2 holds.

\vspace{3mm}\noindent{\bf Lemma 3:~~}%
For any integer $k\geq 1$ and real number $\alpha>0$, there is a
positive integer $C_3$ such that for every integer $n\geq C_3$ and
any positive integer $m<n^\alpha$, we have $(q(m))^k<n$.

\vspace{3mm}\noindent{\bf Proof of Lemma 3:~~}%
Let $r=[\alpha]+1$ be the least integer more than $\alpha$. By Lemma
2, there is a least positive integer $C_4$ such that for every
integer $x\geq C_4$, we have $(q(x))^{kr}<x$. Let $C_5=(C_4+1)^k$.
We will prove that for every integer $n\geq C_5$ and any positive
integer $m<n^\alpha$, we have $(q(m))^k<n$. If $C_4\leq m$, then
$(q(m))^{kr}<m<n^\alpha<n^r$ and $(q(m))^k<n$. If $C_4> m$, then
$(q(m))^k\leq (m+1)^k<(C_4+1)^k=C_5\leq n$. This shows that Lemma 3
holds.

\vspace{3mm}\noindent{\bf  Corollary 1:~~}%
For any given $\varepsilon$ with $0<\varepsilon<0.5$, there is a
positive integer $C_6$ such that for every integer $n\geq C_6$  and
any positive integer $m< n^{2-\varepsilon}$, we have
$2^{\frac{1}{\varepsilon}}(q(m))^{\frac{2-\varepsilon}{\varepsilon}}<n$.

\vspace{3mm}\noindent{\bf  Proof of Theorem  1:~~}%
For any given $\varepsilon$ with $0<\varepsilon<0.5$, there is a
positive integer $C_6$ such that for every prime $p\geq C_6$  and
any positive integer $m< p^{2-\varepsilon}$, we have
$2^{\frac{1}{\varepsilon}}(q(m))^{\frac{2-\varepsilon}{\varepsilon}}<p$
by Corollary 1.

By the prime number theorem in an arithmetic progression, it is easy
to prove that for any prime $p$ with $p\leq\max\{C_1,C_6\}$, ($C_1$
is the positive constant in Conjecture 1), there exists a positive
constant $C_7$ such that for every positive integer $n>C_7$, when
$(p,n)=1$, there exist two distinct odd primes $p_1$ and $p_2$
satisfying $2n\equiv p_1\equiv p_2(\mod p)$ and $p_1,p_2\in
Z_n^*=\{x|1\leq x\leq n, (x,n)=1\}$.

Let $n$ be an integer $>C_7$. Since we assume Conjecture 2, there
exists  $r>1$ such that $p_r<n$, $(p_r,n)=1$ and $2n-p_r$ is coprime
to every $2n-p$ when $p$ ranges the odd primes $\leq n$ and
different from $p_r$. We will show that $2n-p_r$ is prime. If this
is the case, our proof is over, so let us suppose we can write
$2n-p_r =pm$, where $p$ is the least prime factor of $2n-p_r$. Thus,
$2n>p^2$.

We have $p>\max\{C_1,C_6\}$. Indeed, if $p$ is smaller, we can find
two odd primes say $q_1$ and $q_2$, not more than $n$ and prime to
$2n$, such that $2n\equiv q_1\equiv q_2(\mod p)$. At most one of
them, say $q_1$, can be equal to $p_r$. This means that $2n-p_r$ is
not coprime to $2n-q_2$, contrarily to our hypothesis on $p_r$.

Note that $p_r\neq p$ since $(p_r,n)=1$. If $p_r<p$, then
$p+p_r<p^{2-\varepsilon}$ and there is a prime $q$ coprime to
$p+p_r$ and such that
$2^{\frac{1}{\varepsilon}}q^{\frac{2-\varepsilon}{\varepsilon}}<p$
by Corollary 1. Since we suppose that Conjecture 1 holds, hence
there is a prime $x$ such that $x\equiv p+p_r(\mod pq)$ and
$x<(pq)^{2-\varepsilon}<\frac{p^2}{2}<n$. Clearly, $p_r\neq x$. But
$p|(2n-p_r,2n-x)$. It is a contradiction by our assumption on $p_r$.

Hence $p_r>p$. We write $p_r=pl+v$ with $1\leq v<p$. If  $l\geq
p^{1-\varepsilon}$, there is a prime $y$ such that $y\equiv v(\mod
p)$ and $y<p^{2-\varepsilon}<p_r$ (since we suppose Conjecture 1).
But we have also $p|(2n-p_r,2n-y)$, it is contrary to our assumption
on $p_r$ again. So we have $l< p^{1-\varepsilon}$,
$lv<p^{2-\varepsilon}$ and there is a prime $q$ coprime to $lv$ and
such that
$2^{\frac{1}{\varepsilon}}q^{\frac{2-\varepsilon}{\varepsilon}}<p$
by Corollary 1 again. Note that there is a prime $z$ such that
$z\equiv v(\mod pq)$ and $z<(pq)^{2-\varepsilon}<\frac{p^2}{2}<n$
(since we suppose that Conjecture 1 holds). Obviously, we have
$z\neq p_r$ since $(q,l)=1$. But $p|(2n-p_r,2n-z)$. The
contradiction implies that $2n-p_r$ is a prime number. This
completes the proof of Theorem 1.

\vspace{3mm}\noindent{\bf  Remark 1:~~}%
It is interesting that in [11], we proved that if $p(k)< k^2$ and
the necessary condition of Goldbach's conjecture hold, then every
sufficiently large even integer may be written as the sum of a prime
and the product of at most two primes. Namely, our assumptions imply
Chen's theorem [13]. In this paper, we have proved that  if $p(k)\ll
k^{2-\varepsilon}$ and the necessary condition of Goldbach's
conjecture hold, then every sufficiently large even integer may be
written as the sum of two distinct primes. However, it can be
further improved, we think. We hope it can be improved to $p(k)\ll
k^{2+\epsilon}$. Thus, based on work of Chowla [10], one will see
that the Generalized Riemann Hypothesis implies Goldbach's
conjecture. How far "$p(k)\ll k^{2-\epsilon}$" is from "$p(k)\ll
k^{2+\epsilon}$"?

\vspace{3mm} Very naturally, one might ask whether Chowla's
hypothesis is true or not. Of course, due to a limited knowledge of
the author, he can not answer well. However, papers [4] and [6] give
some witnesses. Also based on the structural beauty of Mathematics
itself, the author believe that there is a prime in each row (resp.
column) of the following matrix. This further supports Chowla's
hypothesis.

$$M=(m_{i,j})=\left(
\begin{matrix}
a_1+1\times n, & \cdots,  & a_1+\varphi (n)\times n\\
\cdots, & \cdots, & \cdots \\
a_{\varphi (n)}+1\times n, & \cdots, & a_{\varphi (n)}+\varphi
(n)\times n
\end{matrix} \right),$$ \\where $a_i$ is the $i$-th positive integer
which is coprime to $n$ for $1\leq i\leq \varphi (n)$.

\vspace{3mm}Moreover, for any given integer $n>1$, let $b_i$ be the
$i$-th positive integer which is coprime to $n$ for $1\leq i\leq
\varphi (n)$, where $\varphi (n)$ is Euler totient function, one
might prove that there is a permutation $a_1,...,a_{\varphi (n)}$ of
1 to $\varphi (n)$ such that
$$F_1=\left\{
\begin{array}{c}
f_1(x)=a_1x+b_1\\
...............................\\
f_{\varphi (n)}(x)=a_{\varphi (n)}x +b_{\varphi (n)}\\
\end{array}
\right.$$ is admissible, moreover, $$F_2=\left\{
\begin{array}{c}
f_1(x)=a_1nx+b_1\\
...............................\\
f_{\varphi (n)}(x)=a_{\varphi (n)}nx +b_{\varphi (n)}\\
\end{array}
\right.$$ is admissible, too. Thus, by Dickson's conjecture [15],
$f_1(x),...,f_{\varphi (n)}(x)$ represent simultaneously prime
numbers for infinitely many integers $x$. Therefore, it is very
possible that there is a prime in each row (resp. column) of the
aforementioned matrix. For the definition of `admissible', see [20].

\vspace{3mm} By Chowla's hypothesis, there is a prime in each row of
$M$. By Grimm's conjecture which implies there are two primes
between two square numbers [21], there is a prime in each column of
$M$. From this, we see that many problems in Mathematics are not
isolated again.

\vspace{3mm} On the other hand, we must prove the necessary
condition of Goldbach's conjecture holds without a proviso. This
question looks easy in analytic number theorists' eyes. But, the
author has not been  able to work out a complete proof. We left this
question to the readers who are interested in it. Next, we will try
to consider another interesting problems.

\section{A generalization of the problem of the least prime number in an arithmetic progression}
Clearly, the problem of the least prime number in an arithmetic
progression closely relates to the famous Dirichlet's theorem [14].
In fact, Dirichlet's theorem guarantees us the existence of the
least prime number in  an arithmetic progression.  In 1904, Dickson
[15] generalized Dirichlet's theorem by concerning the simultaneous
values of several linear polynomials,  which implies Green-Tao
theorem [16] (The primes contain arbitrarily long arithmetic
progressions). Unfortunately, Dickson's generalization still is a
conjecture by now.  The author would like to call it Dickson's
conjecture on $N$, where $N$ is the set of all positive integers.

\vspace{3mm}In 2006, Green and Tao [18] considered Dickson's
conjecture in the multivariable case by generalizing
Hardy-Littlewood estimation [22]. The brilliant work Green and Tao
[16] [18] shows that it is possible to generalize Dickson's
conjecture on $N$ to the general case. In [19], the author gave  an
equivalent form of Dickson's Conjecture on $N$ and further
considered Dickson's conjecture on $N^n$. Moreover, in [20], we gave
Dickson's conjecture on $Z^n$ which actually is  an equivalent form
of Green-Tao's conjecture [18], where $Z$ is the set of all
integers.

\vspace{3mm} Well, now, let's assume that Dickson's conjecture on
$N^n$ (or $Z^n$) holds. How to generalize the problem  of the least
prime number in an arithmetic progression? What do  the general
forms of this problem look like?

\vspace{3mm}First,  let's go back to Linnik's theorem [1] again,
which states that $p(k)<k^L$, where $L$ is an absolute constant.
This well-known result can be re-stated as follows: For given
positive integer $k$, there is a positive integer $C_k$ depending on
$k$ such that $p(k)<(C_k)^L$, where $L$ is an absolute constant. Why
do we consider $C_k$? Because for given  $k$ and $C_k$, there are
only finite many integers $l$ satisfying $|l|<C_k$ such that
$f(x)=kx+l$ is admissible. Here $f(x)=kx+l$ represents infinitely
many prime numbers if and only if $f(x)=kx+l$ is admissible. In this
case, we might set $C_k=k$. Based on this simple observation, one
could give a generalization of the problem of the least prime number
in an arithmetic progression as follows:

\vspace{3mm}\noindent{\bf  A generalization of the problem of the
least prime number in an arithmetic progression (A naive approach):~~}%

\vspace{3mm}Let$A=(a_{i,j})=\left(
\begin{matrix}
a_{11}, & \cdots,  & a_{1n}\\
\cdots, & \cdots, & \cdots \\
a_{m1}, & \cdots, & a_{mn}
\end{matrix} \right)$ be an integral matrix in which any two row vectors are not the same
such that for any positive constant $C$, there is an integral point $X=(x_1,...,x_n)\in Z^n$ such
that $F=\left\{
\begin{array}{c}
f_1(X)=a_{11}x_1+...+a_{1n}x_n>C \\
.........................................................\\
f_m(X)=a_{m1}x_1+...+a_{mn}x_n>C\\
\end{array}
\right.$. Then there is a positive integer $C_A$ depending on $A$
such that $p(A)<(C_A)^L$, where $L$ is an absolute constant, and
$p(A)$ be the longest prime vector of $p(A,B)$ for all
$B=(b_1,...,b_m)\in Z^n$ with
$(\sum_{i=1}^{i=m}(b_i)^2)^\frac{1}{2}<C_A$ such that $$G=\left\{
\begin{array}{c}
g_1(X)=a_{11}x_1+...+a_{1n}x_n+b_1 \\
..................................................\\
g_m(X)=a_{m1}x_1+...+a_{mn}x_n+b_m\\
\end{array}
\right.$$ is admissible, and where $p(A,B)$ is the shortest prime
vector (point) represented by $G$. In [17], we have pointed out it
is significative to estimate the upper bound of $p(A,B)$ if $p(A,B)$
exists. What does $C_A$ look like?

\vspace{3mm}\noindent{\bf  Remark 2:~~}%
The condition $F=\left\{
\begin{array}{c}
f_1(X)=a_{11}x_1+...+a_{1n}x_n>C \\
.........................................................\\
f_m(X)=a_{m1}x_1+...+a_{mn}x_n>C\\
\end{array}
\right.$ is necessary because if $G=\left\{
\begin{array}{c}
g_1(X)=a_{11}x_1+...+a_{1n}x_n+b_1 \\
.........................................................\\
g_m(X)=a_{m1}x_1+...+a_{mn}x_n+b_m\\
\end{array}
\right.$ represents infinitely many prime points, then for any
positive constant $C$, there is an integral point
$X=(x_1,...,x_n)\in Z^n$ such that $F=\left\{
\begin{array}{c}
f_1(X)>C \\
...............\\
f_m(X)>C\\
\end{array}
\right.$.

\section{An analogy of  Goldbach's conjecture}
Goldbach's conjecture states that every even integer $2n\geq 8$ is
the sum of two distinct primes. Namely, we have $2n=p+q$, where
$p,q$ are prime with $p<q$. If we look upon $x$ as the value of
number-theoretic function $f(x)$, then when $f(x)=x$, Goldbach's
conjecture can be re-stated as $2f(n)=f(u)+f(v)$ when $n>3$, where
$f(u)=u,f(v)=v$ are prime. Notice that $f(x)=x$ represents
infinitely many prime numbers. More generally, one could expect that
for $\gcd (k,l)=1$, $f(x)=kx+l$ has this property. Namely, for every
sufficiently large integer $n=kw+l$ of the form $f(x)=kx+l$,
$2n=2f(w)=f(u)+f(v)$, where $f(u)=ku+l,f(v)=kv+l$ are prime. This
gives an analogy of  Goldbach's conjecture. For example, let
$f(x)=5x+2$. Then for every $n>9$, $5n+2$ may be written as the sum
of two distinct primes of the form $5x+2$: $2\times 52=104=7+97$,
$2\times 57=114=17+97$, $2\times 62=124=17+107$, $2\times
67=134=7+127$, $144=17+127$, $154=17+137$, $164=7+157$,
$174=17+157$, $184=17+167$, $194=37+157$, $204=7+197$, $214=17+197$,
$224=97+127$, $234=7+227...$

\vspace{3mm} In this section, we give the weakened form of analogy
of Goldbach's conjecture (see Conjecture 3). We further prove that
this weakened form and the weakened form of Chowla's hypothesis
(Conjecture 4) imply the analogy of Goldbach's conjecture.

\vspace{3mm}\noindent{\bf  Conjecture 3:~~}%
Let $k,l$ be given positive integers satisfying $(k,l)=1$ and $1\leq
l<k$. Let $Q_i$ be the $i$th prime of the form $kx+l$. There is a
positive constant $C_8$ such that every integer $n>C_8$, there
exists $r>1$ such that $kn+l>Q_r$, $(kn+l,Q_r)=1$ and $2(kn+l)-Q_r$
is coprime to every $2(kn+l)-Q$ when $Q\leq kn+l$ ranges the primes
of the $kx+l$ and different from $Q_r$.

\vspace{3mm}By Chinese Remainder Theorem, Chowla's hypothesis
implies the following conjecture 4.

\vspace{3mm}\noindent{\bf  Conjecture 4:~~}%
Let $\varepsilon$ with $0<\varepsilon<0.5$ be a real number and
$k,l$ be given positive integers satisfying $(k,l)=1$ and $1\leq
l<k$. There is a positive constant $C_9$, such that for  every
integer $d$ satisfying $d>C_9$, $(d,k)=1$ and any  positive integer
$a$ with $1\leq a<d$ and $(a,d)=1$, there is a prime $q$ such that
$q<(dk)^{2-\epsilon}$, $q\equiv a (\mod d)$ and $q\equiv l (\mod
k)$.

\vspace{3mm}\noindent{\bf  Theorem  2:~~}%
Let $k,l$ be given positive integers satisfying $(k,l)=1$ and $1\leq
l<k$. If Conjecture 3 and Conjecture 4 hold, then for every
sufficiently large integer $n$, $2(kn+l)$ may be written as the sum
of two distinct primes $p,q$ satisfying $p\equiv q\equiv l(\mod k)$.

\vspace{3mm}By the prime number theorem in an arithmetic
progression, more precisely, by the result of Ch. de la
Vall\'{e}e-Poussin [23] which states that $$\sum _{p\equiv l(\mod
k),p\leq x}\log p$$ equals $\frac{x}{\varphi (k)}$ asymptotically,
we have that, for every sufficiently large integer $n$, $\sum
_{i=1}^{i=n+1}\log Q_i$ equals $\frac{Q_{n+1}}{\varphi (k)}$
asymptotically.  It shows immediately that the following Lemma 4
holds.

\vspace{3mm}\noindent{\bf Lemma 4:~~}%
For every integer $h>1$, there is a positive integer $n_h$ such that
$Q_{n+1}^h<Q_1Q_2\cdots Q_n$ for all $n\geq n_h$.

\vspace{3mm}Note that $f(x)=l+kx$ takes infinitely many primes when
$(k,l)=1$. Therefore, for any positive integer $m$, there is a least
prime of the form $l+kx$ which is coprime to $m$. Denote this least
prime by $Q(m)$. By Lemma 4, one can prove the following lemma 5
holds.

\vspace{3mm}\noindent{\bf Lemma 5:~~}%
For every integer $r\geq 1$, there is a positive integer $C_{10}$
depending on $r$, such that for every integer $m\geq C_{10}$, we
have $(Q(m))^r<m$.

\vspace{3mm}\noindent{\bf Lemma 6:~~}%
For any integer $r\geq 1$ and real number $\delta>0$, there is a
positive integer $C_{11}$ depending on $r, \delta$ such that for
every integer $n\geq C_{11}$ and any positive integer $m<n^\delta$,
we have $(Q(m))^r<n$.

\vspace{3mm}\noindent{\bf Proof of Lemma 6:~~}%
Let $e=[\delta]+1$ be the least integer more than $\delta$. By Lemma
5, there is a least positive integer $C_{12}$ such that for every
integer $m\geq C_{12}$, we have $(Q(m))^{re}<m$. Let $Q_g$ be the
least prime of the form $l+kx$ which is larger than $C_{12}$.
Moreover, by Lemma 4, there is a least positive integer $f$ such
that for every integer $h\geq f$, $Q_{h+1}^{re}<Q_1Q_2\cdots Q_h$.
Let $t=\max\{g,f\}$ and $C_{11}=Q_1Q_2\cdots Q_t$. We will prove
that for every integer $n\geq C_{11}$ and any positive integer
$m<n^\delta$, we have $(Q(m))^r<n$.

\vspace{3mm}We write $Q_1Q_2\cdots Q_s\leq n<Q_1Q_2\cdots Q_{s+1}$.
Since $n\geq C_{11}$, hence $s\geq t$. For any positive integer
$m<n^\delta$, if $Q(m)\leq Q_{s+1}$,  then $(Q(m))^{re}\leq
Q_{s+1}^{re}<Q_1Q_2\cdots Q_s$ since $s\geq t\geq f$. So,
$(Q(m))^{re}< n$ and $(Q(m))^r< n$. If $Q(m)> Q_{s+1}$, then $m$ is
divisible by $Q_1Q_2\cdots Q_{s+1}$ and $m\geq Q_1Q_2\cdots
Q_{s+1}$. So, $m>Q_s\geq Q_t\geq Q_g$. Moreover,  $m> C_{12}$ since
$Q_g> C_{12}$. Therefore, $(Q(m))^{re}<m<n^\delta<n^e$ and
$(Q(m))^r<n$. This completes the proof of Lemma 6.

\vspace{3mm}\noindent{\bf  Corollary 2:~~}%
For given $\varepsilon$ in Conjecture 4 and $k$ in Theorem 2, there
is a positive integer $C_{13}$ depending on $\varepsilon,k$ such
that for every prime $p\geq C_{13}$  and any positive integer
$m<k^{3-\varepsilon} p^{2-\varepsilon}$, we have
$2^{\frac{1}{\varepsilon}}(Q(m))^{\frac{2-\varepsilon}{\varepsilon}}k^{\frac{2-\varepsilon}{\varepsilon}}<p$.

\vspace{3mm}\noindent{\bf Lemma 7:~~}%
For given $k$ and $l$ in Theorem 2 and for any odd prime $p$
satisfying $(p,k)=1$ and $p\leq\max\{k+1,C_8,C_{10},C_{13}\}$ (
$C_8,C_{10},C_{13}$ are the aforementioned constants), there exists
a positive constant $C_{14}$ such that for every positive integer
$n>C_{14}$, when $(p,kn+l)=1$, there exist two distinct odd primes
$p_1$ and $p_2$ satisfying $p_1\equiv p_2\equiv l(\mod k)$,
$2(kn+l)\equiv p_1\equiv p_2(\mod p)$ and $p_1,p_2\in
Z_{kn+l}^*=\{x|1\leq x\leq kn+l, (x,kn+l)=1\}$.

\vspace{3mm}\noindent{\bf Proof of Lemma 7:~~}%
Let $D=\max\{k+1,C_8,C_{10},C_{13}\}$ where $C_8,C_{10},C_{13}$ are
the aforementioned constants. For given $k$ and $l$ in Theorem 2,
any odd prime $p$ satisfying $(p,k)=1$ and $p\leq D$, and for any
integer $n$ with $(p,kn+l)=1$, there exists infinitely many primes
$q$ such that $q\equiv l(\mod k)$ and $q\equiv 2(kn+l) (\mod p)$ by
Chinese Remainder Theorem and Dirichlet's prime theorem in an
arithmetic progression. Note that when $n$ ranges positive integers
$\leq p$, $2(kn+l) (\mod p)=\{0,1,2,...,p-1\}$. By the result of Ch.
de la Vall\'{e}e-Poussin  again, for any $r\in \{1,2,...,p-1\}$,
there exists a positive constant $C_{p,r}$ such that for every
positive integer $m$ satisfying $m>C_{p,r}$ and $2(km+l)\equiv r
(\mod p)$, there exist two distinct odd primes $p_1$ and $p_2$
satisfying $p_1\equiv p_2\equiv l(\mod k)$, $2(km+l)\equiv p_1\equiv
p_2(\mod p)$ and $p_1,p_2\in Z_{km+l}^*=\{x|1\leq x\leq km+l,
(x,km+l)=1\}$. Let $C_{14}=\max_{(p,k)=1,p\leq D}\max_{r\in
\{1,2,...,p-1\}}{C_{p,r}}$. It shows that Lemma 7 holds.

\vspace{3mm}\noindent{\bf  Proof of Theorem 2:~~}%
For given $\varepsilon$ in Conjecture 4 and $k,l$ in Theorem 2,
Corollary 2 shows that  there is a positive integer $C_{13}$
depending on $\varepsilon,k$ such that for every prime $p\geq
C_{13}$ and any positive integer $m<k^{3-\varepsilon}
p^{2-\varepsilon}$, we have
$2^{\frac{1}{\varepsilon}}(Q(m))^{\frac{2-\varepsilon}{\varepsilon}}k^{\frac{2-\varepsilon}{\varepsilon}}<p$.

\vspace{3mm}Lemma 7 shows that for any odd prime $p$ satisfying
$(p,k)=1$ and $p\leq\max\{k+1,C_8,C_{10},C_{13}\}$ (
$C_8,C_{10},C_{13}$ are the aforementioned constants), there exists
a positive constant $C_{14}$ such that for every positive integer
$n>C_{14}$, when $(p,kn+l)=1$, there exist two distinct odd prime
$p_1$ and $p_2$ satisfying $p_1\equiv p_2\equiv l(\mod k)$,
$2(kn+l)\equiv p_1\equiv p_2(\mod p)$ and $p_1,p_2\in
Z_{kn+l}^*=\{x|1\leq x\leq kn+l, (x,kn+l)=1\}$.

\vspace{3mm}Let $n$ be an integer $>C_{15}=\max \{C_8,C_{14}\}$.
Since we assume Conjecture 3, there exists $r>1$ such that
$kn+l>Q_r$, $(kn+l,Q_r)=1$ and $2(kn+l)-Q_r$ is coprime to every
$2(kn+l)-Q$ when $Q\leq kn+l$ ranges the primes of the $kx+l$ and
different from $Q_r$. We will show that $2(kn+l)-Q_r$ is prime. If
this is the case, $2(kn+l)-Q_r$ is also a prime of the form $kx+l$
and our proof is over, so let us suppose we can write
$2(kn+l)-Q_r=pm$, where $p$ is the least prime factor of
$2(kn+l)-Q_r$. Thus, $2(kn+l)>p^2$, $(p,kn+l)=1$.

\vspace{3mm}We have $p>\max\{k+1,C_8,C_{10},C_{13}\}$. Indeed, if
$p$ is smaller, we can find two odd primes of the form $kx+l$, say
$q_1$ and $q_2$, not more than $kn+l$ and prime to $2(kn+l)$, such
that $2(kn+l)\equiv q_1\equiv q_2(\mod p)$. At most one of them, say
$q_1$, can be equal to $Q_r$. This means that $2(kn+l)-Q_r$ is not
coprime to $2(kn+l)-q_2$, contrarily to our hypothesis on $Q_r$.

\vspace{3mm}Note that $Q_r\neq p$ since $(kn+l,Q_r)=1$. If $Q_r<p$,
then $(p+Q_r)k<k^{3-\varepsilon}p^{2-\varepsilon}$ and there is a
prime $q$ coprime to $(p+Q_r)k$ and such that
$2^{\frac{1}{\varepsilon}}q^{\frac{2-\varepsilon}{\varepsilon}}k^{\frac{2-\varepsilon}{\varepsilon}}<p$
and $(pqk)^{2-\varepsilon}<\frac{p^2}{2}$ by Corollary 2. Since we
suppose that Conjecture 4 holds, hence there is a prime $A$ of the
form $kx+l$ such that $A\equiv p+Q_r(\mod pq)$ and
$A<(pqk)^{2-\varepsilon}<\frac{p^2}{2}<kn+l$. Clearly, $Q_r\neq A$.
But $p|(2(kn+l)-Q_r,2(kn+l)-A)$. It is a contradiction by our
assumption on $Q_r$.

\vspace{3mm}Hence $Q_r>p$. We write $Q_r=py+z$ with $1\leq z<p$. If
$y\geq p^{1-\varepsilon}k^{2-\varepsilon}$, there is a prime $B$ of
the form $kx+l$ such that $B\equiv z(\mod p)$ and
$B<(pk)^{2-\varepsilon}<Q_r$ (since we suppose Conjecture 4). But we
have also $p|(2(kn+l)-Q_r,2(kn+l)-B)$, it is contrary to our
assumption on $Q_r$ again. So we have $y<
p^{1-\varepsilon}k^{2-\varepsilon}$,
$yzk<p^{2-\varepsilon}k^{3-\varepsilon}$ and there is a prime $q$
coprime to $yzk$ and such that
$2^{\frac{1}{\varepsilon}}q^{\frac{2-\varepsilon}{\varepsilon}}k^{\frac{2-\varepsilon}{\varepsilon}}<p$.
Note that there is a prime $E$ of the form $kx+l$ such that $E\equiv
z(\mod pq)$ and $E<(pqk)^{2-\varepsilon}<\frac{p^2}{2}<kn+l$ (since
we suppose that Conjecture 4 holds). Obviously, we have $E\neq Q_r$
since $(q,y)=1$. But $p|(2(kn+l)-Q_r,2(kn+l)-E)$. The contradiction
implies that $2(kn+l)-Q_r$ is a prime number. Therefore, when
$n>C_{15}$, $2(kn+l)$ may be written as the sum of two distinct
primes $p,q$ satisfying $p\equiv q\equiv l(\mod k)$ assuming
Conjectures 3 and 4. This completes the proof of Theorem 2.

\vspace{3mm}\noindent{\bf  Remark 3:~~}%
Based on Euclid's algorithm, in [27], we find a special sequence
which is called W sequence. By studying W sequences in the case of
non-consecutive positive integers, we give Conjectures 2 and 3.
Conjecture 4 can be generalized: Let $\varepsilon$ with
$0<\varepsilon<0.5$ be a real number and $k_i,l_i$ be given positive
integers satisfying $(k_i,l_i)=1$ and $1\leq l_i<k_i$ for
$i=1,...,n$, where $(k_i,k_j)=1$ for $1\leq i\neq j\leq n$. There is
a positive constant $C_{16}$, such that for every integer $d$
satisfying $d>C_{16}$, $(d,k_1...k_n)=1$ and any positive integer
$a$ with $1\leq a<d$ and $(a,d)=1$, there is a prime $q$ such that
$q<(dk_1...k_n)^{2-\epsilon}$, $q\equiv a (\mod d)$ and $q\equiv l_i
(\mod k_i)$.  This can be deduced by  Chinese Remainder Theorem and
Chowla's hypothesis.

\section{ A generalization of analogy of  Goldbach's conjecture (A naive approach)}
\vspace{3mm}It is known that $f(x)=x$ on $Z$ is the simplest
polynomial which represents infinitely many primes. By  Dirichlet's
famous theorem, for any positive integer $l,k$ with $(l,k)=1$,
$f(x)=l+kx$ is a simpler polynomial which also represents infinitely
many primes. If we view  $f(x)=l+kx$ as an analogy of $f(x)=x$,
Theorem  2 shows that it is possible to  give an analogy of
Goldbach's conjecture. Lev Landau said tastily: `Why add prime
numbers? Prime numbers are made to be multiplied, not added.' This
time, we are afraid of `Prime numbers might be made to be added.' If
in the higher-dimension case, we have a similar the problem  of the
least prime number in an arithmetic progression, maybe, there is
also a similar Goldbach's conjecture. We will try to consider this
problem in this section.  Very naturally, We would like to consider
a general problem: might prime points be made to be added? In order
to clearly explain this problem, firstly, let's do an interesting
thing as follows:

\vspace{3mm} Based on the point of view that a number is a map, we
view an integer $x$ as the simplest polynomial map on $Z$: $f(x)=x$
from $Z$ to $Z$.  Notice that such a map takes infinitely many prime
numbers. More generally, let's consider the map $F: Z^n\rightarrow
Z^m$ for all integral points $x=(x_1,...,x_n)\in Z^n$,
$F(x)=(f_1(x),...,f_m(x))$ for distinct polynomials $f_1,...,f_m\in
Z[x_1,...,x_n]$, where $m,n\in N$. In this case, we call $F$ a
polynomial map on $Z^n$. We say that these multivariable integral
polynomials $f_1(x),...,f_m(x)$ on $Z^n$  represent simultaneously
prime numbers for infinitely many integral points $x$, if for any
$1\leq i\leq m$, $f_i(x)$ itself can represent prime numbers for
infinitely many integral points $x$, moreover, there is an infinite
sequence of integral points $(x_{11},...,x_{n1})$, ...,
$(x_{1i},...,x_{ni})$, ... such that for any positive integer $r$,
$f_1(x_{1r},...,x_{nr})$,..., $f_m(x_{1r},...,x_{nr})$ represent
simultaneously prime numbers, and for any $i\neq j$,
$f_1(x_{1i},...,x_{ni})\neq f_1(x_{1j},...,x_{nj})$, ...,
$f_m(x_{1i},...,x_{ni})\neq f_m(x_{1j},...,x_{nj})$ hold
simultaneously. In this case, we also say that the polynomial map
$F$ on $Z^n$  represents infinitely many prime points. Such a
polynomial map $F$ is called a prime map. In short, a prime map is a
polynomial map on $Z^n$ which represents infinitely many prime
points.  For instance: $F=f(x)=x$, $F=f(x)=ax+b$ with $\gcd
(a,b)=1$, $F=f(x,y)=x^2+y^2$, $F=f(x,y)=x^2+y^2+1$,
$F=f(x,y)=x^3+2y^3$, $F=f(x,y)=x^2+y^4$,
$F=(f_1(x,y)=x,f_2(x,y)=x^2+y^2)$, $F=f(x,y,z,w)=x^2+y^2+z^2+w^2$
and so on are all prime maps. This gives a generalization of
$f(x)=x$ on $Z$. Due to the fact that $g(x)=ax+b$ with $\gcd(a,b)=1$
is the unique known prime map on $Z$, we want to know more
properties about the arithmetic progressions. By the analogy of
Goldbach's conjecture, we further hope to find more interesting
analogies between Integers and Arithmetic progressions. For example,
for every sufficiently large integer $n$, if $g(x)>n$, then there is
a prime of the form $g(x)$ in the interval $(g(x),2g(x))$, which can
be viewed as the analogy of Bertrand-Chebyshev theorem, especially,
there is a prime of the form $3k+1$ in the interval $(3x+1,2(3x+1))$
for each positive integer $x$. These problems we will study in other
papers.

\vspace{3mm}In [20], we find an interesting property of prime maps
and generalize the analogy of Chinese Remainder Theorem as follows:
Let $F=(f_1,...,f_m)$ be a prime map. If $\gcd (a_i,a_j)=1$ for
$1\leq i\neq j\leq k$, and there exist integral point $x^{(j)}\in
Z^n$ such that $F(x^{(j)})$ is in $(Z_{a_j}^*\setminus \{1\})^m$ for
$1\leq j\leq k$, then there exists an integral point $z$ such that
$F(z)$ is in $(Z_{a_1...a_k}^*\setminus \{1\})^m$.

\vspace{3mm}Note that the prime map
$F=(f_1(x,y)=x,f_2(x,y)=x^2+y^2)$ implies that for any $1\leq m\leq
n$, there is a prime map $F(x)=(f_1(x),...,f_m(x))$ for distinct
polynomials $f_1,...,f_m\in Z[x_1,...,x_n]$. However, when $m>n$, we
do not know whether there are always such prime maps. Especially,
when $n=1,m>1$ and $f_i$ is linear, it is a famous open problem
(Dickson's conjecture). Anyway, one might expect that prime maps
have many fascinating properties like integers. We expect that prime
points have some interesting properties like prime numbers. The
author wishes that in the higher-dimension case, we have a similar
Prime Number Theorem.

\vspace{3mm} We call a prime map $F(x)=(f_1(x),...,f_m(x))$ on $Z^n$
is standard if $F(1,...,1)$ is a prime point (vector). For example,
$f(x)=x+1$ is a standard prime map. $F=f(x)=ax+b$  a standard prime
map if and only if $a+b$ is a prime number. Bertrand-Chebyshev
theorem implies that for any positive integer $a>1$, there is a
positive integer $b<a$ such that $ax+b$ is a standard prime map.
Clearly, a prime maps can be reduced to a standard prime map. Let
$F=(f_1,...,f_m)$ on $Z^n$ be a standard prime map. Then for every
sufficiently large integer $r$, if there is an integral point $x$
such that each coordinate of $F(x)=\alpha$ is greater than $r$, then
$2\alpha=\beta+\gamma$, where $\beta,\gamma$ are distinct prime
points represented by $F$. This explains the aforementioned problem.
From this, one will see that this problem and the prime map are
equivalent. Particularly, Goldbach's conjecture and the infinitude
of primes are equivalent although without any proof. This perhaps is
another property of primes maps. But, this problem is the author's
naive viewpoint. The author also finds several propositions which
are equivalent to the infinitude of prime numbers by considering
prime maps, see [Appendix]. For this reason, we revisit Euclid's
Number Theory and focus on the essence of integers. G\"{o}del's
incompleteness theorem [25] states that all consistent axiomatic
formulations of number theory include undecidable propositions.
Along this research line, we do not know whether one will meet those
undecidable propositions in Number Theory. We hope that people
further consider them.

\section{Acknowledgements}
I am very thankful to  Professor Heath-Brown for his comments
improving the presentation of the paper, and also to my supervisor
Professor Xiaoyun Wang for her  help. Thank Professor Xianmeng Meng
for her suggestions. Thank the key lab of cryptography technology
and information security in Shandong University and the Institute
for Advanced Study in Tsinghua University, for providing me with
excellent conditions. This work was partially supported by the
National Basic Research Program (973) of China (No. 2007CB807902)
and the Natural Science Foundation of Shandong Province (No.
Y2008G23).

\section* {Appendix: Euclid's Number Theory Revisited}

\textbf{Remark:} This appendix is  self-contained.

\vspace{3mm} From Euclid's famous \emph{Elements} [24], (Proposition
20, Book ¢ù), we see that Euclid (300 B.C.) proved that $f(x)=x$
represents infinitely many prime numbers.

\vspace{3mm} It is difficult to image what would happen if there was
only finite many prime numbers: many theorems and conjectures do not
hold any more. For example,  Bertrand-Chebyshev theorem, Dirichlet's
theorem, Prime number theorem, The  fundamental theorem of
arithmetic, Chinese Remainder Theorem, Goldbach's conjecture,
Landau's problems and so on are not true if there is only finite
many prime numbers. Therefore, Hardy said: `Euclid's theorem which
states that the number of primes is infinite is vital for the whole
structure of arithmetic. The primes are the raw material out of
which we have to build arithmetic, and Euclid's theorem assures us
that we have plenty of material for the task'. For this reason, in
this section, we would like to revisit Euclid's Number Theory and
give some equivalent propositions of Euclid's second theorem.

\vspace{3mm} From Book 7, 8 and 9 of his \emph{Elements}, we see
that Euclid had established elementarily Theory of Divisibility and
the greatest common divisor.  Euclid began his number-theoretical
work by giving some definitions and his algorithm (the Euclidean
algorithm) (See [24]: Book 7, Propositions 1 and 2) as follows:

......

11. A prime number is that which is measured by a unit alone.

12. Numbers prime to one another are those which are measured by a
unit alone as a common measure.

13. A composite number is that which is measured by some number.

......

Proposition 1 (Book 7): Two unequal numbers being set out, and the
less being continually subtracted in turn from the greater, if the
number which is left never measures the one before it until a unit
is left, the original numbers will be prime to one another.

Proposition 2 (Book 7): Given two numbers not prime to one another,
to find their greatest common measure.

......

Proposition 31 (Book 7): Any composite number is measured by some
prime number.

......

Proposition 20 (Book 9): Prime numbers are more than any assigned
multitude of prime numbers. Namely, there are infinitely many
primes.

......

\vspace{3mm} Now, let's go back to Euclid's proof for the infinitude
of prime numbers: Supposed that there are only finitely many primes,
say $k$ of them, which denoted by $p_1,...,p_k$. Consider the number
$E=1+ \prod_{i=1}^{i=k}p_i$. If $E$ is prime, it leads to the
contradiction since $E\neq p_i$ for any $1\leq i\leq k$. If $E$ is
not prime, $E$ has a prime divisor $p$ by Proposition 31 (Book 7).
But $p\neq p_i$ for any $1\leq i\leq k$. Otherwise, $p$ divides
$\prod_{i=1}^{i=k}p_i$. Since it also divides
$1+\prod_{i=1}^{i=k}p_i$, it will divide the difference or unity,
which is impossible.

\vspace{3mm} In his proof, we see that Euclid used Proposition 31
(Book 7). Of course, he also used a unexpressed axiom which states
that if $A$ divides $B$, and also divides $C$, $A$ will divide the
difference between $B$ and $C$.

\vspace{3mm}Well, let's look at the proof of Proposition 31 (Book
7): Let $A$ be a composite number. By the definition, there must be
a number $B$ ($1<B<A$) which divides $A$. If $B$ is prime, then
Proposition 31 holds. If  $B$ is not prime, there must be a number
$C$ ($1<C<B$) which divides $B$. If $C$ is prime, then Proposition
31 holds since $C$ also divides $A$. If $C$ is not prime, by
repeating this process, in  finite many steps, there must be a prime
which divides $A$ and Proposition 31 holds. From this proof, we see
that Euclid used a unexpressed axiom which states that if $A$
divides $B$, and $B$ divides $C$, then $A$ divides $C$. In his book
[24], Thomas Little Heath had noted that Euclid used the
aforementioned axioms. We would be quite surprised if he did use
these axioms because  on one hand, Proposition 31 (Book 7) and
Proposition 20 (Book 9) can be deduced early by definitions, on the
other hand, we expect him to make use of his algorithm which is his
first number-theoretical proposition in his \emph{Elements}. Then,
let's try to supplement this work.

\vspace{3mm} Now, let's use these axioms again, Euclid's definitions
on a prime number and a composite number, and his algorithm
(Propositions 1, 2) to prove Euclid's second theorem and some
equivalent propositions of the infinitude of prime numbers.

\vspace{3mm}\noindent {\bf Theorem 1:~~}%
Any composite number is divided by some prime number.

\vspace{3mm}\noindent {\bf Proof:~~}%
By the definition of a composite number and the axiom which that if
$A$ divides $B$, and $B$ divides $C$, then $A$ divides $C$, it is
easy to prove that Theorem 1 is true.

\vspace{3mm}\noindent {\bf Theorem 2:~~}%
For any positive integer $a$, $a$ is co-prime to $a+1$.

\vspace{3mm}\noindent {\bf Proof:~~}%
By Euclid's algorithm (Proposition 1, Book 7), it shows immediately
that Theorem 2 holds.

\vspace{3mm}\noindent {\bf Corollary 1:~~}%
For any positive integer $a$, there is a positive integer $b$ such
that $b>1$ and $b$ is co-prime to $a$.

\vspace{3mm}\noindent {\bf Proof:~~}%
Let $b=a+1$. By Theorem 2, it is easy to prove Corollary 1 holds.

\vspace{3mm}\noindent {\bf Corollary 2:~~}%
For any positive integer $a$, there is a least integer $b$ such that
$b>1$ and $b$ is co-prime to $a$.

\vspace{3mm}\noindent {\bf Proof:~~}%
By Corollary 1 and the axiom which states there is a least element
in any non-empty subset of natural numbers, Corollary 2 holds.

\vspace{3mm}\noindent {\bf Theorem 3:~~}%
For any positive integer $a$, let $b$ be the least integer $b$ such
that $b>1$ and $b$ is co-prime to $a$. Then $b$ is prime.

\vspace{3mm}\noindent {\bf Proof:~~}%
By Corollary 2, we get that for any positive integer $a$, there is a
least integer $b$ such that $b>1$ and $b$ is co-prime to $a$. If $b$
is not prime, by Theorem 1, $b$ is divided by some prime number $p$.
Of course,  $p$ is co-prime to $a$ and $p<b$. But $b$ is the least.
The contradiction shows that $b$ is prime and Theorem 3 holds.

\vspace{3mm}\noindent {\bf Corollary 3:~~}%
2 and 3 are all  prime numbers.

\vspace{3mm}\noindent {\bf Proof:~~}%
We do not want to factor 2 or 3 but prove directly Corollary 3
holds. By Theorem 2, we know that $1$ is co-prime to $2$. Note that
2 is the least integer such that $2>1$ and $2$ is co-prime to $1$.
Let $a=1$. By Theorem 3, we deduce that 2 is prime. Similarly, one
can prove that 3 also is prime.

\vspace{3mm} Corollary 3 gives us a method for generating whole
prime numbers: Let $p_i$ be the $i$-th prime. By Corollary 3,
$p_1=2$, $p_2=3$. $p_{n+1}$ is the least prime which is co-prime to
$\prod_{i=1}^{i=n}p_i$.

\vspace{3mm}\noindent {\bf Corollary 4:~~}%
There are infinitely many prime numbers.

\vspace{3mm}\noindent {\bf Proof:~~}%
The existence of prime numbers is very clear. for example, 2 is a
prime number by Corollary 3. Supposed that there are only finitely
many prime numbers, say $k$ of them, which denoted by $p_1,...,p_k$.
Let $a=\prod_{i=1}^{i=k}p_i$. By Theorem 3, let $b$ be the least
integer $b$ such that $b>1$ and $b$ is co-prime to $a$. Then $b$ is
prime. Of course, $b\neq p_i$ for any $1\leq i\leq k$. The
contradiction shows that Corollary 4 is true.

\vspace{3mm} From Corollary 3, we see that Propositions 1 and 31
(Book 7) in Euclid's Number Theory implies the infinitude of prime
numbers. Next, we will give some equivalent propositions that there
are infinitely many prime numbers. The author wonders why this
occurs.

\vspace{3mm}\noindent {\bf Theorem 4:~~}%
There are infinitely many prime numbers if and only if for any
positive integer $a$, there is a positive integer $b$ such that
$b>1$ and $b$ is co-prime to $a$.

\vspace{3mm}\noindent {\bf Proof:~~}%
If there are infinitely many prime numbers, then for any positive
integer $a$, there must be a prime $p$ which is greater than $a$.
Let $b=p$ and the necessity holds obviously. On the other hand, if
for any positive integer $a$, there is a positive integer $b$ such
that $b>1$ and $b$ is co-prime to $a$, then there must be a least
integer $c$ such that $c>1$ and $c$ is co-prime to $a$. By Theorem
3, $c$ is prime. Thus the existence of prime numbers has been
proved. Supposed that there are only finitely many prime numbers,
say $k$ of them, which denoted by $p_1,...,p_k$. Let
$d=\prod_{i=1}^{i=k}p_i$. By Theorem 3 again, let $e$ be the least
integer such that $e>1$ and $e$ is co-prime to $d$. Then $e$ is
prime. Of course, $e\neq p_i$ for any $1\leq i\leq k$.  The
contradiction shows that the sufficiency is true. Therefore, Theorem
4 holds.

\vspace{3mm} From Theorem 4, we see that the polynomial $f(x)=x$
takes infinitely many prime numbers if and only if it is admissible.

\vspace{3mm}\noindent {\bf Lemma 1:~~}%
Euclid's algorithm, Division algorithm and Bezout's equation are all
equivalent.

\vspace{3mm}\noindent {\bf Proof:~~}%
See [26].

\vspace{3mm} Since we aforehand assume Euclid's algorithm, hence, by
Lemma 1, we can logically deduce many number theoretical results in
any number theoretical textbooks. Especially, we get the following
theorems 5, 6 and 7.

\vspace{3mm}\noindent {\bf Theorem 5:~~}%
For any positive integer $a$, there is a positive integer $b$ such
that $b>1$ and $b$ is co-prime to $a$ if and only if there is a
positive integer $c$ such that for any positive integer $m>c$, there
is a positive integer $k$ such that $1<k<m$ and $k$ is co-prime to
$m$.

\vspace{3mm}\noindent {\bf Proof:~~}%
First, we prove that the latter implies the former. When $a>c$,
since for any positive integer $m>c$, there is a positive integer
$k$ such that $1<k<m$ and $k$ is co-prime to $m$, hence there is a
positive integer $b$ such that $b>1$ and $b$ is co-prime to $a$.
When $1<a\leq c$, clearly, there is a positive integer $r$ such that
$a^r>c$. Thus, there is a positive integer $b$ such that $b>1$ and
$b$ is co-prime to $a^r$. Of course, $b$ is co-prime to $a$, too.
When $a=1$, we can choose $b=2$.

\vspace{3mm} Next, we will prove that the former implies the latter.
One might believe that for any positive integer $m>2$, there is a
positive integer $k=m-1$ such that $1<k<m$ and $k$ is co-prime to
$m$ by Theorem 2. Thus, it seems that the former is not related to
the latter. However, we do not do so. We will strictly prove that
for any positive integer $m\geq 15$, there is a positive integer $k$
such that $1<k<m$ and $k$ is co-prime to $m$. Clearly, if 3 (resp.
5) is co-prime to $m$, we choose $k=3$ (resp. $k=5$). So, when
$m\geq 15$, we only consider the case that $m$ is divisible by 15.
We write $m=15t$ with $t\geq1$. If $t$ is not divisible by 2, we can
choose $k=2$. Well, now we assume that $t$ is divisible by 2. We
write $t=3^ed$ with $\gcd (3,d)=1$. Since $2|t$, hence $d>1$. Note
that there is a positive integer $r>1$ which is co-prime to $3t$
because we assume that for any positive integer $a$, there is a
positive integer $b$ such that $b>1$ and $b$ is co-prime to $a$. By
the linear congruence theorem, there is a positive integer $h$ with
$0\leq h<d\leq t$ such that $2+3h\equiv r (\mod d)$. Notice that
either $2+3h$ or $2+3h+3t$ is co-prime to $m=15t$, moreover,
$1<2+3h<2+3h+3t<15t=m$. Let $c=15$. This shows that Theorem 5 holds.

\vspace{3mm}\noindent {\bf Theorem 6:~~}%
There are infinitely many prime numbers if and only if there is a
positive constant $c$ such that for any positive integer $a>c$,
there is a positive integer $b$ such that $1<b<a$ and $b$ is
co-prime to $a$.

\vspace{3mm}\noindent {\bf Proof:~~}%
By Theorems 4 and 5, it immediately shows that Theorem 6 is true.

\vspace{3mm} From Theorem 6, we see that the polynomial $f(x)=x$
takes infinitely many primes if and only if it is strongly
admissible.

\vspace{3mm}\noindent {\bf Corollary 5:~~}%
There are infinitely many prime numbers if and only if for any
positive integer $a>2$, there is a positive integer $b$ such that
$1<b<a$ and $b$ is co-prime to $a$.

\vspace{3mm}\noindent {\bf Proof:~~}%
By Theorem 5, the infinitude of prime numbers implies for any
positive integer $a>14$, there is a positive integer $b$ such that
$1<b<a$ and $b$ is co-prime to $a$. Further, one can directly test
that it  is also true when $2< a< 15$. So Corollary 5 holds.

\vspace{3mm}\noindent {\bf Theorem 7:~~}%
Euclid's second theorem and the analogy of Chinese Remainder Theorem
(which states that if  there exist a positive integer $a$ such that
$1<a$ is in $Z_n^*$ and  a positive integer  $b$ such that $1<b$ is
in $Z_m^*$, then there exists a positive integer $c$ such that $1<c$
is in $Z_{mn}^*$ when $\gcd (m,n)=1$) are equivalent.

\vspace{3mm}\noindent {\bf Proof:~~}%
First, we prove that  the analogy of Chinese Remainder Theorem
implies Euclid's second theorem. By Corollary 3, we know that 2 and
3 are all prime numbers. Supposed that there are only finitely many
prime numbers, say $k$ of them, which denoted by
$p_1=2,p_2=3,...,p_k$. Let $d=\prod_{i=2}^{i=k}p_i$. Clearly, 2 is
in $Z_d^*$. 3 is in $Z_4^*$. Notice that $\gcd (d,4)=1$. By our
assumption, there exists a least positive integer $c$ such that
$1<c$ is in $Z_{4d}^*$. By Theorem 3, $c$ is prime. Of course,
$c\neq p_i$ for any $1\leq i\leq k$.  The contradiction shows that
the analogy of Chinese Remainder Theorem implies Euclid's second
theorem.

\vspace{3mm}Secondly, we prove that Euclid's second theorem implies
the analogy of Chinese Remainder Theorem. By Corollary 5, we only
need to prove that for any positive integer $d>2$ there is a
positive integer $k$ such that $1<k<d$ and $\gcd (k,d)=1$ implies
the analogy of Chinese Remainder Theorem. In fact, if  there exist a
positive integer $a$ such that $1<a$ is in $Z_n^*$ and  a positive
integer $b$ such that $1<b$ is in $Z_m^*$, then $m\geq 3,n\geq 3$.
Consequently $mn\geq 9>2$. By our assumption that Euclid's second
theorem holds, equivalently, when $mn>2$, there is a positive
integer $c$ such that $1<c<mn$ and $\gcd (c,mn)=1$, we deduce that
the analogy of Chinese Remainder Theorem holds. This completes the
proof of Theorem 7.

\vspace{3mm}  One might find more equivalent propositions.  By the
aforementioned discussion, we believe that one of substantive
characteristics of the set of all integers is that it contains
infinitely many prime numbers. Therefore, it should be reasonable
that we generalize Integers to Prime maps.  Based on such a belief,
we revisited Euclid's Number Theory and added this appendix.

\clearpage
\end{document}